\documentclass[amssymb,amscd,12pt]{article}
\usepackage{epsfig}
\usepackage{amsmath}
\usepackage{amsthm}
\usepackage{amssymb}

\usepackage{amscd}
\usepackage{amsfonts}  
\usepackage{makeidx}
\usepackage{hyperref}
\usepackage[all]{xy}
\usepackage{graphicx}    
\usepackage{multicol}    

\DeclareMathAlphabet\EuR{U}{eur}{m}{n}
\SetMathAlphabet\EuR{bold}{U}{eur}{b}{n}

\newcounter{commentcounter}
\newcommand{\comment}[1]                      
{\stepcounter{commentcounter}
{\bf Comment \arabic{commentcounter}}: {\ttfamily #1} }

\newtheorem{theorem}{Theorem}

\newtheorem{proposition}[theorem]{Proposition}

{\catcode`@=11\global\let\c@equation=\c@theorem}

\title {Involutions on 3-Manifolds and Self-dual, Binary Codes
}
\author{Matthias Kreck and Volker Puppe}
\begin{document}
\maketitle
\noindent



\vspace{1cm}

{\footnotesize {\bf Abstract} We study a correspondence between orientation
reversing involutions on compact 3-manifolds with only isolated
fixed points and binary, self-dual codes. We show in particular
that every such code can be obtained from such an involution. We
further relate doubly even codes to $Pin^-$-structures and
$Spin$-manifolds.}

\section{Introduction}

Binary, self-dual codes play an important role in coding theory and have been studied extensively, see [RS] for a comprehensive survey and
literature. A connection to involutions on 3-manifolds was made in [Pu]. It is shown there that an involution $\tau: M \to M$ on a
$3$-dimensional, closed manifold $M$ with ''maximal number'' of
isolated, fixed points (i.e. with only isolated fixed points such
that the number of fixed points $k: = |M^{\tau}|$ equals $dim_
 {\mathbb  Z/_2} (\oplus_i H^i (M; \mathbb  Z/_2)) $ determines a
binary self-dual code of length $k$. In turn this code determines
the cohomology algebra $H^* (M; \mathbb  Z/_2)$ and the equivariant
cohomology $H^*_G (M; \mathbb  Z/_2)$, where the action of $G: =
\mathbb Z/_2$ is given by the involution. In fact, the code
corresponds to the inclusion $H^*_G (M; \mathbb  Z/_2)
\longrightarrow H^*_G (M^G; \mathbb  Z/_2) \cong (\mathbb  Z/_2
[t])^k$, and an equivalence of codes given by permuting the
coordinates, corresponds to an equivalence of inclusions, given by
an automorphism of the algebra
$ (\mathbb  Z/_2 [t])^k$. \\

In Section 2 we generalize these results by considerung involutions on 3-manifolds which have a finite number of fixed points which need not
be maximal. The code corresponding to such an involution is described in two ways, firstly using equivariant cohomology as in [Pu], secondly
using the ordinary homology (with $\mathbb Z/_2$ coefficients) of the complement(of a neighbourhood) of the fixed point in the orbit space. \\

In Section 3 we show that every binary, self-dual code can be obtained from an involution on a 3-manifold with a finite number of
fixed points,
in fact, using surgery we get that the manifold can be chosen so that the number of fixed points is maximal.\\

In Section 4 we relate doubly even codes to $Spin$-manifolds. We define the concept of  a $Spin$-involution and show that $Spin$-involutions
give doubly even codes. Finally we show that each doubly even code comes from a $3$-manifold with $Spin$-involution.

\section{Self-dual codes from involutions on 3-manifolds}

 Let $\tau: M \to M$ be an involution on a closed $3$-manifold $M$ with finitely
many fixed points $x_1, \dots, x_k$. By Smith theory $k \le
dim_{\mathbb  Z/_2}(\oplus_i H^i (M; \mathbb  Z/_2))$. The famous
localization Theorem for equivariant cohomology gives in this
context that the map
$$
H^*_G(M) \stackrel{i^*_G} {\longrightarrow} H^*_G (M^G) \cong
(\mathbb  Z/_2 [t])^k
$$
(here and in the following we always take
coefficients in $\mathbb  Z/_2$), induced by the inclusion $M^G
\stackrel{i}{ \longrightarrow} M$, becomes an isomorphism after
inverting
the powers of $t \in \mathbb  Z/_2 [t] = H^* (BG)$. \\

The kernel of $i^*_G$ is the torsion submodule $T \subset H^*_G
(M)$; in particular $i^*_G$ is injective if and only if $H^*_G
(M)$ is a free $\mathbb  Z/_2[t]$- module. As in the ''maximal
number of fixed point'' case, the inclusion $H^*_G (M)/T \stackrel
{\bar{i}^*_G}{\longrightarrow} (\mathbb  Z/_2 [t])^k$ determines a
self-dual code $C$, and in turn $\bar{i}^*_G$ is determined by
this code. In more detail $C$ is the image of $\stackrel
{1}{\oplus}_{i=0} H^i_G (M)$ in $\mathbb  Z/_2^k$ under the
evaluation map, putting $t=1$ (cf. [Pu]). This is the same as the
image of $H^1_G (M)$ in $H^1_G(M^G) \cong \mathbb  Z/_2^k$.  The
code $C$ is self-dual. This is shown in [Pu] for the ''maximal
case'' (i.e. $T = 0$), and can be seen generally as follows (cf.
[AP], Exercise (1.15)): \\
Let $\rho : H^*_G (M) \to H^*(M)$ be the restriction to the fibre
in the Borel construction $M \to M \times_G EG \to BG$.  The map
$\rho$ fits into the long exact Gysin sequence of the covering $M
\simeq M \times EG \to M \times_G EG$, i.e.
$$\dots \to H^{i-1}_G (M) \stackrel {\cup t} {\to} H_G^i (M) \stackrel{\rho}
{\to} H^i(M)  {\to} H_G^i(M) \stackrel {\cup t} {\to} H_G^{i+1} (M) \to \dots$$
is exact and hence so is
$$
0 \to H_G^* (M) \otimes_{\mathbb Z/_2 [t]} \mathbb  Z/_2 \to H^*
(M)
 \to Tor ^ {\mathbb  Z/_2 [t] } (H_G^* (M), \mathbb  Z/_2) \to 0.
$$
Splitting $H^*_G (M) \cong T \oplus F$ into a direct sum of the
torsion submodule $T$ and a free complement $F$, one sees that
$dim_{\mathbb  Z/_2} (T \otimes_{\mathbb  Z/_2 [t]} \mathbb  Z/_2) =
dim_{\mathbb  Z/_2} Tor ^ {\mathbb  Z/_2 [t] } (H_G^* (M), \mathbb
\mathbb Z/_2)$. We claim that with respect to the Poincar\'e duality
pairing in $H^* (M)$,
the orthogonal complement of  $ \rho (H^*_G (M))$ is $\rho (T)$.\\

Because of the dimension equality above it suffices to show that
$<x,y> = \sigma (x \cup y) = 0$ for $x \in \rho (H^*_G(M))$ and $y
\in \rho(T)$, where $\sigma : H^*(M) \to \mathbb  Z/_2$ is the
orientation of $M$. If $\tilde {x}$ and $\tilde{y}$ are liftings
of $x$ and $y$ with respect to $\rho$, then $\tilde {x} \cup
\tilde{y}Ê\in T$. Therefore $\tilde {x} \cup \tilde{y}$ is mapped
to zero under the equivariant orientation $\tilde{\sigma} : H^*_G
(M) \to \mathbb  Z/_2[t]$ (cf. e.g. [AP], Chap. 1). Hence
$\tilde{\sigma} (\tilde {x} \cup \tilde{y}) = 0$
and thus $ \sigma (x \cup y) =0$.\\

We get that the graded algebra $§ H^*_G (M)/_§ T$  is a
subquotient of $H^*(M)$, which fulfills Poincar\'e duality, and
$H^*_G (M) /T$ as a $\mathbb  Z/_2 [t]$-module is isomorphic to $(
\rho (H^*_G(M))/ \rho (T)) \otimes_{\mathbb  Z/_2} \mathbb
Z/_2[t]$. We therefore have, analogous to the case $T = 0$, that
$$
H^*_G (M)/T
\stackrel{{\bar{i}}^*_G}{\longrightarrow} (\mathbb  Z/_2 [t])^k
$$
determines  a binary, self-dual code and, in turn, is determined
by this code. Note, though, that we only get the quotient algebra
$H^*_G(M)/T$ (and $\rho (H^*_G(M))/ \rho(T) \cong (H^*_G(M)/T)
\otimes_{\mathbb  Z/_2 [t]} \mathbb  Z/_2$ ) from the code which
means that in case $T \not= 0$ the algebras $H^*_G (M)$ and $H^*
(M)$ are not
completely determined by the code.\\

In view of the construction below, we describe the code coming
from an involution $\tau : M \to M$ on a closed 3-manifold $M$
with isolated fixed points $x_1, \dots, x_k$ in a second way.\\

Let $W: = (M \setminus +_k \stackrel{\circ} {D^3})/\tau$, where
$D^3_i$ are equivariant discs around $x_i$. We consider the
Mayer-Vietoris sequence in equivariant cohomology for $M = (M
\setminus +_k \stackrel{\circ}{D^3}) \cup (+_k D^3): \dots \to
H^0_G ( +_k S_1^2) \to H^1_G (M) \to H^1_G (M \setminus +_k
\stackrel{\circ}{D^3}) \oplus H^1_G (+_k D^3) \to H^1_G (+_k S^2)
\to \dots$ \\since the equivariant cohomology of a free $G$-space
is the non-equivariant cohomology of the orbit space, one gets the
exact sequence
$$
\dots \to H^0 (+_k \mathbb R P^2) \to H^1_G(M) \to H^1 (W) \oplus
(H^1_G(+_k D^3) \to H^1(+_k \mathbb R P^2) \to \dots.
$$
It is easy  to see that the map $ (H^1_G(+_k D^3) \to H^1(+_k
\mathbb R P^2)$ is an isomorphism. (For one disk $D^3$ one has
that $S^2 \stackrel{(id, i)}{\longrightarrow} S^2 \times
S^{\infty}  \stackrel{(j, id)}{\longrightarrow} D^3 \times
S^\infty \stackrel {p}{\to} S^\infty$ is the canonical inclusion,
and hence $H^1_G (D^3) \cong H^1 (\mathbb RP^{\infty}) {\to}
H^1(\mathbb RP^2) = H^1_G(S^2)$. It therefore follows from the
Mayer-Vietoris sequence that $Im (H^1_G (M) \to H^1_G (M^G)) = Im
(H^1 (W) \to H^1(+_k \mathbb RP^2))$. Dually to taking $Im (H^1
(W) \to H^1 (+_k \mathbb RP^2))$, we can take $Ker (H_1 (+_k
\mathbb RP^2) \to H_1 (W))$. Since the kernel of the map on the
middle homology of the boundary of a compact manifold to the inner
is a self annihilating subspace (with respect to the intersection
form) of half rank and the intersection form on $H_1 (+_k \mathbb
RP^2) = (\mathbb  Z/_2)^k$ is the standard ''Euclidean'' form, one
gets in another way that the code
is self-dual.\\

We summarize the above considerations somewhat vaguely as

\begin{theorem} Every involution with only isolated fixed points on a compact 3-manifold determines a binary, self-dual code.
\end{theorem}

\section{All self-dual codes come from 3-manifolds}

\begin{proposition}
Every binary, self dual code can be obtained from an involution on an orientable 3-manifold.
\end{proposition}

\noindent
{\bf Proof:}
Let $k = 2r$. Let $C \subset {\mathbb  Z/_2}^{k}$ be a self dual code. We choose a map $f: +_k \mathbb {RP}^2 \to (\mathbb  {RP}^\infty)^r$ such
that the sequence (with $\mathbb{Z}_2$ coefficients)
$$
0 \to C \to H_1(+_k \mathbb {RP}^2)Ê\stackrel{f_*}{ \to}  H_1(( \mathbb {RP}^\infty)^r) \to 0
$$
is exact. \\

Next we note that the first Stiefel-Whitney class $w_1(+_k \mathbb {RP}^2)$ is in the image $f^*$. The reason is that the diagonal
element $\Delta$ is in the code ($\Delta$ is dual to  $w_1(+_k \mathbb {RP}^2)$) and
so $<w_1(+_k \mathbb {RP}^2),x> = <\Delta  , x>  =0$ for
all $x\in C$. This implies that there is a real line bundle  $L$  over $(\mathbb  {RP}^\infty)^r$ pulling back to
the non-trivial line bundle
over each copy of $\mathbb {RP}^2$. Thus $(+_k \mathbb {RP}^2,f)$ is an orientable singular manifold, where the orientation is
twisted by
the line bundle $L$. This means that the bundle $\nu (+_k \mathbb {RP}^2) - f^*(L)$ is orientable.\\

After choosing an orientation the pair $(+_k \mathbb{R}P^2,f)$ represents an element in the bordism
group $(\Omega_2(({\mathbb{R}P^{\infty}})^r;L)$ of singular manifold with orientation (twisted by $L$). We claim that
this element is trivial. \\

The Atiyah-Hirzebruch spectral sequence implies that
$$
\Omega_2 (( \mathbb {RP}^\infty)^{r}; L)
\stackrel{\cong}{\longrightarrow} H_2(( \mathbb {RP}^\infty)^{r};
\mathbb Z_t).
$$
$$
[F,h] \longrightarrow h_* ([F])
$$
Here, $\mathbb Z_t$ stands for twisted homology, where the coefficient system is given by the representation
$$
\pi_1 ((\mathbb {RP}^\infty)^{r} \to \pi_1 ( \mathbb {RP}^\infty) \stackrel{\cong}{\longrightarrow} (\pm 1) = Aut (\mathbb Z)
$$
and the map is induced by the classifying map of $L$.
We note that
$$
H_2  ((\mathbb {RP}^\infty)^{r}; \mathbb Z_t) \to H_2 (( \mathbb {RP}^\infty)^{r}; \mathbb Z/_2)
$$
is injective. The reason is that  $H_2  ((\mathbb {RP}^\infty)^{r}; \mathbb Z_t)$ consists only of elements of order $2$ (and $0$).\\

Thus it is enough to control the image of  the fundamental class in $H_2((\mathbb{RP}^\infty)^{r}; \mathbb Z/_2)$. The vanishing is equivalent
to
$$
h^* x \cup h^* y = 0
$$
for all $x, \; y \in H^1 ((\mathbb{RP}^\infty)^{r}; \mathbb
Z/_2)$.  This follows since by construction the intersection form
vanishes on $Ker (H_1 (+_k \mathbb{RP}^2; \mathbb Z/_2) \to H_1
((\mathbb{RP}^\infty)^{r}; \mathbb Z/_2)$, which under the
isometry between $H_1 (+_k \mathbb{RP}^2; \mathbb Z/_2) \cong H^1
(+_k \mathbb{RP}^2; \mathbb Z/_2)$ corresponds to the
image of
 $H^1 ((\mathbb{RP}^\infty)^{r}; \mathbb Z/_2) \to H^1(+_k \mathbb{RP}^2; \mathbb Z/_2)$. \\

Summarizing the information so far, we have shown that
$$
[+_k \mathbb {RP}^2, f] = 0 \in \; \Omega_2((\mathbb
{RP}^\infty)^{r}; L).
$$\\

Let $h : W \to (\mathbb {RP}^\infty)^{r}$ be a zero bordism. We claim that the kernel of the map induced by the inclusion
$$
H_1 (+_k \mathbb {RP}^2; \mathbb Z/_2) \to H_1 (W; \mathbb Z/_2)
$$
is our code $C$. Since $h|_{+_k \mathbb {RP}^2} = f$, we conclude
that $Ker (H_1 (+_k \mathbb {RP}^2; \mathbb Z/_2) \to H_1 (W;
\mathbb Z/_2)$
is contained in $C$. But this kernel has dimension $r = dim \; C$ implying the statement. \\

Finally we consider the classifying map $g $ of $L$ and the composition
$$
g h : W \to \mathbb {RP}^\infty,
$$
 to construct the induced 2-fold covering $\hat{W}$ over $W$.  Since $W$ is oriented  (twisted by $ L$),  $\hat{W}$ is an orientable manifold.
 The boundary of $W$ is $+_k S^2$ and the restriction of the deck transformation to the boundary is $-id$ on each summand. Thus we obtain an
 involution $\tau$ on
$$
M:=  \hat{W} \cup  +_k  D^3
$$
which on $\hat{W}$ is the deck transformation and on each $D^3$ is $-id$. By construction the code associated to this 3-manifold $M$
and $\tau$ is the given code finishing the argument.\\
{\bf q.e.d.}\\

The above construction depends on the choice of the zero cobordism
$h:W \to (\mathbb{R}P^{\infty})^r$, and it is not clear whether
one obtains a manifold $M$ with involution, which has maximal
number of isolated fixed points. We will show that one can change
$W$ by surgery to reduce the cohomology of $M$, and obtain a pair
$(M, \tau)$ with maximal number of isolated fixed points.  By
Smith theory the maximality condition is equivalent to the
injectivity of $H^*_G (M; \mathbb{Z}_2) \to H^*_G (M^G;
\mathbb{Z}_2) \cong (\mathbb{Z}_2 [t])^k$, resp. the surjectivity
of $H^G_* (M^G; \mathbb{Z}_2) \to H^G_*(M; \mathbb{Z}_2)$. In our
case $M= \hat{W} \cup (+_k D^3)$. The equivariant Mayer-Vietoris
sequence (with $\mathbb{Z}_2$ coefficients) gives:
$$
\dots \to H^G_*(+_k S^2) \to H^G_* (\hat{W}) \oplus H^G_*(+_k D^3) \to H^G_*(M) \to H^G_{*-1}(+_k S^2) \to \dots
$$
One has $H^G_*(+_k S^2) \cong (H_*(+_k \mathbb{R}P^2)$ and
$H^G_*(\hat{W}) \cong H_* (W)$ since the actions on $+_k S^2$ and
$\hat{W}$ are free; and $H^G_* (+_k D^3) \cong H^G_* (M^G)$. The
inclusion $S^2 \subset D^3$ induces the inclusion $H_*
(\mathbb{R}P^2) \to H_*(\mathbb{R}P^{\infty})$. Hence the map
$H^G_*(M^G) \to H^G_*(M)$ is surjective if and only if $H_i (+_k
\mathbb{R}P^2) = H_i (\partial W) \to H_i (W) $ is surjective for
$i=1,2$. But the long exact sequence of the Poincar\'e pair $(W,
\partial W)$ gives that surjectivity for $i=1$ already implies
surjectivity for $i=2$. To verify the maximality condition it
therefore suffices to show that $H_1(\partial W) \to H_1(W)$ is
surjective. We want to arrive at this condition by applying
surgery to $W$ (if necessary). Assume that $H_1(\partial W) \to
H_1(W)$ is not surjective. We consider the following diagram
$$
\begin{array}{cccccc }
 \to  & H_1(\partial W) & \stackrel{i_1}{\longrightarrow} & H_1(W) & \to & H_1(W, \partial W)  \\
      &   f_1 \searrow & & \swarrow h_1  & & \\
      &  &  H_1((\mathbb{R}P^{\infty})^r) & & & \\
\end{array}
$$
We already know that $i_1$ and $f_1$ have the same kernel, namely
the code $C$, and $f_1$ is surjective by construction. Hence $i_1$
is surjective if and only if $h_1$ is injective. Assume that there
exists an $a \in H_1(W), a \not= 0$, with $h_1(a)=0$. Since $W$ is
orientable (twisted by $L$), the normal bundle of an embedded
circle representing $a$ is trivial. Performing surgery with
respect to $\alpha$ kills the class $a$ and its dual with respect
to the intersection pairing. The map $h_1:H_1(W) \to
H_1((\mathbb{R}P^{\infty})^r)$ factors through the quotient
$H_1(W)/_{<a>}$. Hence we can find a map $h':W' \to
(\mathbb{R}P^{\infty})^r$ of our new manifold $W'$, which
restricts to $f$ on the
 boundary $\partial W'= \partial W = +_k \mathbb{R}P^2$. Iterating the process (if necessary) gives the following result.

\begin{theorem}
Every binary, self-dual code can be obtained from an involution on an orientable 3-manifold with maximal number of isolated fixed points.
\end{theorem}

\bigskip
\bigskip
\bigskip
\bigskip



\section{Spin-structures and doubly even codes}

Let $M$ be an oriented closed $3$-manifold with involution $\tau$
having exactly $k$ isolated fixed points . We will construct from
these data a $4$-manifold by starting with $M \times S^1$ and
dividing by the involution, which on $M$ is $\tau$ and on $S^1$ is
complex conjugation. This is a manifold with $2k$ isolated
singularities, where $k$ is the number of fixed points of $\tau$.
All fixed points singularities are cones over $\mathbb {RP}^3$,
which are the links of the singularities. Since the involution on
$M \times S^1$ is orientation preserving, the orientation on $M
\times S^1$ induces an orientation on the quotient (after removing
the fixed points), which in turn gives an orientation on each link
$\mathbb {RP}^3$. Now we remove open cones around the
singularities and replace them by the disc bundle of the complex
line bundle over $\mathbb {CP}^1$ with Chern class $-2$. The
reason for choosing this sign of the Chern class (and not $+2$) is
that the induced orientation on $\mathbb {RP}^3$ above is the
opposite of this orientation (we will discuss this in more detail
in the proof of the following result). This implies that the
orientations fit together  and so the result is an oriented
$4$-manifold
denoted $N(M,\tau)$. We say that $\tau$  is a $Spin$-involution if $N(M,\tau)$ admits a Spin structure compatible with the given
orientation.\\

The construction of $N(M,\tau)$ is well known in the case of the $3$-torus $T^3$ with $\tau$ complex conjugation.
Then $N(T^3,\tau)$ is the $K3$-surface, which has a Spin structure.

\begin{theorem} Let $M$ be a closed oriented $3$-manifold with involution $\tau$ with finitely many fixed points.
If $\tau$ is a $Spin$-involution, then the code $C(M,\tau)$ is doubly even.
\end{theorem}

\noindent
{\bf Proof:} We assume that the reader is familiar with $Pin^-$-structures [KT]. We recall that a $Pin^-$-structure on a smooth
manifold $M$ is a $Spin$-structure on  $TM \oplus Det (TM)$. Here we note that $TM \oplus Det (TM)$ has a natural orientation,
which we assume to be compatible with the $Spin$-structure. Thus the $Pin^-$-structures are classified by $H^1(M;\mathbb Z/_2)$. \\

A $Pin^-$-structure on a surface $F$ determines a quadratic
refinement $q: H^1(F;\mathbb Z/_2) \to \mathbb Z /4$, such that
$q(x+y) = q(x) + q(y) + 2<x,y>$, where $<x,y>$ is the intersection
form. The two $Pin^-$-structures on $\mathbb {RP}^2$ are
distinguished by the quadratic form, which can take the values
$\pm 1$. For all this see [KT]. If $W$ is a $3$-dimensional
$Pin^-$-manifold with boundary $F$, then on the image of
$H^1(W;\mathbb Z/_2) \to H^1(\partial W;\mathbb Z/_2)$ the
intersection form and the quadratic refinement vanish. This
follows from [KT] as explained in [T].
\\

Now suppose that the disjoint union of $k$ copies of $\mathbb
{RP}^2$
 is the boundary of a
 $Pin^-$-manifold $W$
and the induced $Pin^-$-structure is equal on all components of
the boundary, then if $x\in H^1(W;\mathbb Z/_2)$, and $ y = i^*(x)
= (y_1,\dots,y_{k})$ we conclude that $\sum_i y_i$ = 0 mod $4$.
Thus we are finished if the condition that $\tau$ is a
$Spin$-involution implies that $(M \setminus \ +_k {B^3})/\tau$
has a $Pin^-$-structure, which on all boundary components
is the same. Here ${B^3}_i$ is an open ball around the $i$-th fixed point. \\

To see this we first note that a $Pin^-$-structure on $\mathbb
{RP}^2$ is the same as a $Spin$-structure on the total space of
$T\mathbb {RP}^2 \oplus Det\,T\mathbb {RP}^2$. Since $Det T\mathbb
{RP}^2$ is the normal bundle of $\mathbb {RP}^2$ in $\mathbb
{RP}^3$, we can via a tubular neighbourhood identify  $T\mathbb
{RP}^2 \oplus Det\,T\mathbb {RP}^2$ with an open subset of
$\mathbb {RP}^3$, which is homotopy equivalent to $\mathbb {RP}^3
- pt$. Thus a $Pin^-$-structure on $\mathbb {RP}^2$ determines a
$Spin$-structure on $\mathbb {RP}^3$ and vice versa. In particular
this means that the $Pin^-$-structure on $\mathbb {RP}^2$
determines an orientation on $\mathbb {RP}^3$. We note that
$\mathbb {RP}^3$ is the total space of the complex line bundle
over $\mathbb {CP}^1$ with first Chern class $2$. Using the
complex orientation on $\mathbb {CP}^1$ and on the complex line
bundle we obtain an orientation on $\mathbb {RP}^3$. It is not
difficult to show that this orientation agrees with the
orientation coming from the $Pin^-$-structure on $\mathbb {RP}^2$
(one only has to compare the orientations at one point). Thus, if
this is the orientation on a component of the boundary of some
$4$-manifold $V$, then we obtain an oriented manifold by gluing
the disk bundle of the complex line bundle over $\mathbb{CP}^1$
with Chern class $-2$ (this induces the negative orientation
compared to the orientation above, and so the orientations fit
together). The key observation for our proof is, that since
this disc bundle is simply connected there is a unique $Spin$-structure on it.  \\

Now we consider $M \times S^1$ with the involution given by $\tau$
and complex conjugation $c$. Each fixed point of $M$ and each
fixed point of $S^1$ gives a fixed point of $M \times S^1$ and for
each fixed point the link of the corresponding singularity in $M
\times S^1/(\tau \times c)$ is $\mathbb {RP}^3$ containing the
link of the corresponding singularity in $M/\tau$. Thus a
$Pin^-$-structure on each link in $M /\tau$ determines a
$Spin$-structure of the two (for each fixed point of $S^1$)
corresponding links  in $M \times S^1/(\tau \times c)$) and vice
versa.  If $N(M,\tau)$ has a $Spin$-structure, this  is the same
on each disk
 bundle of the complex line bundle with Chern class $-2$ over $\mathbb {CP}^1$, since
 there is a unique $Spin$-structure with the given orientation. Thus the restriction of the $Pin^-$-structure to each
 link $\mathbb {RP}^2$ is the same.
 As explained above this implies the theorem.\\
{\bf q.e.d.}\\

Next we prove that for each doubly even self dual code $C$ there is a $3$-manifold $M$ with $Spin$-involution $\tau$ such
that the corresponding code is $C$.

\begin{proposition} Let $C$ be a doubly even self dual code. Then there is a $3$-manifold $M$ with $Spin$-involution $\tau$ whose
code is $C$.
\end{proposition}

\noindent {\bf Proof:} We proceed as in the proof of Theorem 3 and
use the notation from there. Now we consider $+_k\mathbb {RP}^2$
as a $Pin^-$-manifold, where all copies are equipped with the same
$Pin^-$-structure, which, if we pass to the corresponding
$Spin$-structure on $\mathbb {RP}^3$ can be extended to the disc
bundle of the complex line
 bundle with first Chern class $-2$ over $\mathbb {CP}^1$. Together with the map $f$ we obtain an element
 of $\Omega_2^{Pin^-}( \mathbb {RP}^\infty)^r$. We compute this bordism group with the Atiyah-Hirzebruch spectral sequence.
 We use from [KT] that $\Omega_0^{Pin^-} = \mathbb Z$, $\Omega_1^{Pin^-} = \mathbb Z/_2$ and $\Omega_2^{Pin^-} = \mathbb Z/8$, where the
 latter group is generated by $\mathbb {RP}^2$ with any $Pin^-$-structure.\\

The Atiyah-Hirzebruch spectral sequence computing
$\Omega_2^{Pin^-}( (\mathbb {RP}^\infty)^r)$ has the entries:
$$\Omega_2^{Pin^-},$$
$$ H_1(( \mathbb {RP}^\infty)^r;\mathbb Z/_2),$$
$$H_2 ((\mathbb {RP}^\infty)^r;\mathbb Z).$$
The component in the first entry is given by the $k$-fold sum of
$\mathbb {RP}^2$ with the given $Pin^-$-structure, which is zero
if and
only if $k = 0 \,mod\,8$. But this is the case for doubly even self dual codes. \\

The last entry is as in the case of oriented bordism (twisted by $L$) detected by the image of the fundamental class with
coefficients in $\mathbb Z/_2$, which - as shown before - vanishes if the code is self dual.\\

The second entry is a bit delicate. We only have to detect the
corresponding entry for $\Omega_2^{Pin^-}( \mathbb {RP}^\infty)$,
since we can project to the different components. Then the
corresponding entry is in $\mathbb Z/_2$. By the fact that the
bordism group is a module over $\Omega_*^{Pin^-}$ we see that the
non-trivial element is represented by $(S^1 \times \eta,  ip_1)$,
where $\eta $ is $S^1$ with the non-trivial $Pin^-$-structure
(which for $1$-manifolds is the same as a $Spin$-structure) and
$i$ is the inclusion $S^1 \to \mathbb {RP}^\infty$. We are free to
choose a $Pin^-$-structure  on the first factor. If we choose the
$Spin$-structure again to be the non-trivial one, we see that the
induced $2$-fold cover is $\eta \times \eta$, which is the
non-trivial element in $\Omega _2 ^{Spin}$. We note that whatever
$Pin^-$-structure we choose on the first factor, we can change it,
if necessary, to the non-trivial one, by modifying it with the
non-trivial element in the image of $H^1(\mathbb
{RP}^\infty;\mathbb Z/_2)$. The upshot of these considerations is
that we can detect the second term in the Atiyah-Hirzebruch
spectral sequence of an element $[N,g] \in \Omega
_2^{Pin^-}(\mathbb {RP}^\infty)$ whose underlying $Pin_-$-bordism
class is zero and whose fundamental class maps to zero  by the
following criterion: It is zero, if and only if  for all
modifications of the $Pin^-$-structure by elements in $H^1(\mathbb
{RP}^\infty;\mathbb Z/_2)$ the induced $2$-fold covering is zero
bordant.  Applying this to the case, where $N  = +_{8l}
\mathbb{RP}^2$ we note that the induced covering is an $S^2$'s
over each summand which maps
 non-trivial to $\mathbb {RP}^\infty$ (which is zero bordant) and that
 it is $\mathbb{RP}^2+ \mathbb{RP}^2$ for each summand which
 maps trivially. But since $\mathbb{RP}^2$ is a generator
 of $\Omega _2^{Pin_-} \cong \mathbb Z/8$, this implies that if the number
 of summands which are  mapped trivial is $0  mod  4$, then
 the bordism class is trivial. Returning to the situation given by our code
 we see that if the code is doubly even, this criterion applies.\\

Thus we have shown that for doubly even codes the bordism class
vanishes in $\Omega _2((\mathbb {RP}^\infty)^r)$, and as in the
proof of Theorem 2 we construct a $3$-manifold $M$ with involution
$\tau$ giving the code. Since the $Pin^-$-structure is the same
for all copies of $\mathbb {RP}^2$, we obtain a $Spin$-involution.
Namely the $4$-manifold we construct is the blow up of a
$Spin$-manifold obtained by replacing the open cones over the
individual $\mathbb {RP}^3$'s by the disc bundle of the complex
line bundle with Chern class $-2$ over $\mathbb {CP}^1$. After
perhaps changing the orientation before the gluing, the resulting
manifold is oriented.
Since the $Spin$-structure on all $\mathbb {RP}^3$'s extend to this disc bundle, the manifold is a $Spin$-manifold.\\
{\bf q.e.d.}\\

As before one can apply surgery, this time taking into account the
$Pin^-$-structure to get the following result.\\

\begin{theorem}
Every binary, doubly even, self-dual code can be obtained from a
$Spin$-involution with maximal number of isolated fixed points on
an orientable 3-manifold.
\end{theorem}


\bigskip
\bigskip

The first author was partially supported by a DFG grant. He also would like to thank the Mathematics Department of Columbia
University and the Courant Institute for support and a stimulating atmosphere while part of the project was carried out.

\bigskip
\bigskip

\noindent
[AP] Allday, C. and  Puppe, V.: Cohomological Methods in transformation Groups. cambridge Studies in Advanced Math. 32, Cambridge University Press, Cambridge, 1993\\
\noindent
[KT] R. C. Kirby and L. R. Taylor.: Pin structures on low-dimensional manifolds.
Geometry of Low-dimensional Manifolds, 2, London Math. Soc. Lecture Note Ser., 151, Cambridge Univ. Press (1990) 177-242.  \\
\noindent
[Pu] Puppe, V.: Group Actions and Codes, Canad. J. Math., Vol 53 (1), 2001, 212 - 224.\\
\noindent
[RS] Rains, E. M. and Sloane, N. J. A.: Self-dual codes, in Handbook of Coding Theory, Pless, V. S. and Huffman, W. C., Eds., 
Elsevier, Amsterdam, 1998, 177-294.\\
\noindent
[T] L.R.Taylor: $Pin^-$-structures, e-mail (2006)\\  

\bigskip
\bigskip

\end {document}